%% file: banachtarski.tex
\documentclass[11pt]{article}
 \input{mathdefs}

\begin{document}
\title{The Banach--Tarski paradox \\ 
or \\
What mathematics and religion have in common}
\author{\it Volker Runde}
\date{}
\maketitle
\begin{quote} \begin{small}
As he went ashore he saw a great throng; and he had compassion on them, and healed their sick. When it was evening, the disciples came to him and said:
``This is a lonely place, and the day is now over; send the crowds away to go into the villages and buy food for themselves.'' Jesus said: ``They need not go away;
you give them something to eat.'' They said to him: ``We have only five loaves here and two fish.'' And he said: ``Bring them here to me.'' Then he ordered
the crowds to sit down on the grass; and taking the five loaves and the two fish he looked up to heaven, and blessed, and broke and gave the loaves to the disciples,
and the disciples gave them to the crowds. And they all ate and were satisfied. And they took up twelve baskets full of the broken pieces left over. And those
who ate were about five thousand men, besides women and children.
\par\hfill {\sf Mt 14:14--21}
\end{small} \end{quote}
\par
Why does an article that is supposed to be about mathematics start with the feeding of the five thousand?
\par
In the nineteen twenties, two Polish mathematicians --- Stefan Banach and Alfred Tarski --- proved a mathematical theorem which sounds a lot like the feeding of the 
five thousand. In their honor, it is called the {\it Banach--Tarski paradox\/}\footnote{The theorem is proved in the article: {\sc S.\ Banach} and {\sc A.\ Tarski}, Sur 
la d\'ecomposition des ensembles de points en parts respectivement congruents. {\it Fund.\ Math.\/} {\bf 6\/}
(1924), 244--277. Look it up if you think I'm making things up!}. Consequences of the Banach--Tarski paradox are, for example:
\begin{quote}
An orange can be chopped into a finite number of chunks, and these chunks can then be put together again to yield {\it two\/} oranges, each of which {\it has the same size as the one
that just went into pieces\/}.
\end{quote}
Another, even more bizarre consequence is:
\begin{quote}
A pea can be split into a finite number of pieces, and these pieces can then be reassembled to yield a solid ball whose diameter is larger than the distance of the Earth
from the sun.
\end{quote}
More generally, whenever you have a three-dimensional body (with a few strings attached), you can obtain any other such body by breaking the first one 
into pieces and the reassembling the parts. To turn five loaves and two fish into enough food to feed a crowd of more than five thousand, then just appears to be a minor exercise.
\par
If you have read this far, your attitude will presumably one of the following:
\begin{itemize}
\item You are a fundamentalist Christian and are delighted to find that mathematics lends support to your beliefs.
\item You are a fundamentalist Christian and are infuriated of how mathematicians make a mockery of your beliefs.
\item You are not a fundamentalist Christian, but your belief in the absolute truth of mathematical theorems is so strong that it makes you swallow
the Banach--Tarski paradox.
\item You are a person who puts common sense above everthing else, so that you neither take the feeding of the five thousand nor the Banach--Tarski paradox at face value.
\end{itemize}
If you fall into the first three categories, there is probably little incentive for you to read any further. Otherwise, I guess, your attitude is best described as
follows: You may believe in story of the feeding of the five thousand, but not take it literally, and if you hear of a mathematical theorem whose consequences are obviously nonsense,
you tend to think that the theorem is wrong.
\par
Take an orange, a sharp knife, and a chopping block. Then chop the orange into pieces, and try to form two globes of approximately the same size out of the orange chunks. 
If the chunks are small enough, each of these two globes will bear reasonable resemblance to a ball, but, of course, each of them has a volume which is only about half of that
of the orange you started with. Perhaps you just didn't chop up the orange in the right way? Give it another try. The result will be the same. You can try your luck on 
hundreds of oranges: You will produce tons of orange pulp, but no corroboration of the Banach--Tarski paradox. Doesn't this show that the Banach--Tarski paradox is wrong? 
\par
The Banach--Tarski paradox is a so-called existence theorem: There is a way of splitting up a pea such the pieces can be reassembled into, say, a life-sized statue of
Stefan Banach. The fact that you haven't succeeded in finding such a way doesn't mean that it doesn't exist: You just might not have found it yet. Let me clarify this with an
example from elementary arithmetic. A positive integer $p$ is called {\it prime\/} if $1$ and $p$ itself are its only divisors; for example, $2$, $3$, and $23$ are prime, whereas
$4 = 2 \cdot 2$ and $243 = 3 \cdot 81$ aren't. Already the ancient Greeks knew that every positive integer has a {\it prime factorization\/}: If $n$ is a positive 
integer, then there are prime numbers $p_1, \ldots, p_k$ such that $n = p_1 \cdot \cdots \cdot p_k$. For small $n$, such a prime factorization is easy to find: $6 = 2 \cdot 3$,
$243 = 2 \cdot 3 \cdot 3 \cdot 3 \cdot 3$, and $6785 = 5 \cdot 23 \cdot 59$, for example. There is essentially only one way of finding a prime factorization: trying. Already finding 
the prime factorization of $6785$ --- armed only with pencil and paper --- would have taken you some time. And now think of a large number, I mean, {\it really large\/}:
\[
  7380563434803675764348389657688547618099805.
\]
This is a perfectly nice positive integer, and the theorem tells you that it has a prime factorization, but --- please! --- don't waste hours, days, or even years of your life 
trying to find it. You might think: What were computers invented for? It is easy to write a little program that produces the prime factorization of an arbitrary
positive integer (and it may even produce one of $7380563434803675764348389657688547618099805$ in a reasonable amount of time). The avarage time, however, it takes such a program
to find the prime factorization of an integer $n$ goes up dramatically as $n$ gets large: For sufficiently large $n$, even the fastest super-computer available today would
--- in avarage --- require more time to find the prime factorization of $n$ {\it than the universe already exists\/}. So, although a prime factorization of a positive integer
always exists, it may be impossibly hard to find. In fact, this is a good thing: It is at the heart of the public key codes that make credit card transactions
on the internet safe, for example. Now, think again of the Banach--Tarski paradox: Just because you couldn't put it to work in your kitchen (just as you couldn't find
the prime factorization of some large integer), this doesn't mean that the theorem is false (or that this particular integer doesn't have a prime factorization).
\par
Let's thus try to refute the Banach--Tarski paradox with the only tool that works in mathematics: pure thought. What makes the Banach--Tarski paradox defy common sense is
that, apparently, the volume of something increases out of nowhere. You certainly know a number of formulae to calculate the volumes of certain particular three-dimensional bodies: 
If $C$ is a cube whose edges have length $l$, then its volume $V(C)$ is $l^3$; if $B$ is a ball with radius $r$, then its volume $V(B)$ is $\frac{4}{3}\pi r^3$.
But what's the volume of an arbitrary three-dimensional body? No matter how the volume of a concrete body is calculated, the following are certainly true about the volumes
of arbitrary, three-dimensional bodies:
\begin{items}
\item if the body $\tilde{B}$ is obtained from the body $B$ simply by moving $B$ in three-dimensional space, then $V(\tilde{B}) = V(B)$;
\item if $B_1, \ldots, B_n$ are bodies in three-dimensional space, then the volume of their union is less than or equal to the sum of their volumes, i.e.\
\[
  V(B_1 \cup \cdots \cup B_n) \leq V(B_1) + \cdots + V(B_n).
\]
\item if $B_1, \ldots, B_n$ are bodies in three-dimensional space such that any two of them have no point in common, then the volume of their union is the even equal to the 
sum of their volumes, i.e.\
\[
  V(B_1 \cup \cdots \cup B_n) = V(B_1) + \cdots + V(B_n).
\]
\end{items}
So, let $B$ be an arbitrary three-dimensional body, and let $B_1, \ldots, B_n$ be subsets of $B$ such that any two of them have no point in common and
$B = B_1 \cup \cdots \cup B_n$. Now, move each $B_j$ in three-dimensional space, and obtain $\tilde{B}_1, \ldots, \tilde{B}_n$. Finally, put the $\tilde{B}_j$ together and
obtain another body $\tilde{B} = \tilde{B}_1 \cup \cdots \cup \tilde{B}_n$. Then we have for the volumes of $B$ and $\tilde{B}$:
\begin{eqnarray*}
  V(B) & = & V(B_1 \cup \cdots \cup B_n) \\
  & = & V(B_1) + \cdots + V(B_n), \qquad\text{by (iii)}, \\
  & = & V(\tilde{B}_1) + \cdots + V(\tilde{B}_n), \qquad\text{by (i)}, \\
  & \geq & V(\tilde{B}_1 \cup \cdots \cup \tilde{B}_n),\qquad\text{by (ii)},\\
  & = & V(\tilde{B}).
\end{eqnarray*}
This means that the volume of $\tilde{B}$ must be less than or equal to the volume of $B$ --- it can't be larger. Banach and Tarski were wrong! Really?
\par
Our refutation of Banach--Tarski seems to be picture perfect: All we needed were three very basic properties of the volume of three-dimensional bodies. But was this really all?
Behind our argument, there was a hidden assumption: Every three-dimensional body has a volume. If we give up that assumption, our argument suddenly collapses: If only one of
the bodies $B_j$ has no volume, our whole chain of (in)equalities makes no longer sense. But why shouldn't every three-dimensional body have
a volume? Isn't that obvious? What is indeed true is that every orange chunk you can possibly produce on your chopping block has a volume. For this reason, you will never be able
to use the Banach--Tarski paradox to reduce your food bill. A consequence of the Banach--Tarski paradox is therefore: There is a way of chopping up an orange, so that you can form, say, 
a gigantic pumpkin out of the pieces --- but you will never be able to do that yourself using a knife. What kind of twisted logic can make anybody put up with that?
\par
Perhaps, you are more willing to put up with the {\it  axiom of choice\/}:
\begin{quote}
If you have a family of non-empty sets $S$, then there is a way to choose an element $x$ from each set $S$ in that family.
\end{quote}
That sounds plausible, doesn't it? Just think of a finite number of non-empty sets $S_1, \ldots, S_n$: Pick $x_1$ from $S_1$, then proceed to $S_2$, and finally take
$x_n$ from $S_n$. What does the axiom of choice have to do with the Banach--Tarski paradox? As it turns out, a whole lot: {\it If\/} the axiom of choice is true, then 
the Banach--Tarski paradox can be derived from it, and, in particular, there must be three-dimensional bodies without volume. So, the answer to the question of whether the 
Banach--Tarski paradox is true depends on whether the axiom of choice is true.
\par
Certainly, the axiom of choice works for a finite number of non-empty sets $S_1, \ldots, S_n$. Now think of a sequence $S_1, S_2, \ldots$ of non-empty sets. Again, pick
$x_1$ from $S_1$, then $x_2$ from $S_2$, and just continue. You'll never come to an end, but eventually you'll produce some element $x_n$ from each $S_n$. So, the axiom of choice
is true in this case, too. But what if we have a truly arbitrary family of sets? What if we have to deal with the family of all non-empty subsets of the real line? It can be shown that this 
family of sets can't be written as a sequence of sets. How do we pick a real number from each set? There is no algorithm that enables us to pick one element from one set, then to proceed to 
another and to eventuall pick an element out of each of them. Nevertheless, the axiom of choice still seems plausible: Each set $S$ in our family is non-empty and therefore contains some 
element $x$ --- why shouldn't there be a way of choosing a particular element from each such set?
\par
I hope, I have convinced you that the axiom of choice is plausible. On the other hand, it implies strange phenomena like the Banach--Tarski paradox. If it's true we have to
put up with the mysterious duplication of oranges. If it's false, then: Why? Please, don't try to prove or to refute the axiom of choice. You won't succeed either way.
The axiom of choice is beyond proof or refutation: We can {\it suppose\/} that it's true, or we may {\it suppose\/} that it's false. In other words: We have to {\it believe\/}
in it or leave it. Most mathematicians these days are believers in the axiom of choice for a simple reason: With the axiom of choice they can prove useful theorems, most of which are 
much less baffling than the Banach--Tarski paradox.
\par
Are you disappointed? Instead of elevating the feeding of five thousand from a matter of belief to a consequence of a bullet-proof mathematical theorem, the Banach--Tarski paradox
demands that you accept another article of faith --- the axiom of choice --- before you can take the theorem for granted. After all, the Banach--Tarski paradox is not that much
removed from the feeding of the five thousand\ldots
\vfill
\begin{tabbing}
{\it Address\/}: \= Department of Mathematical Sciences \\
\> University of Alberta \\
\> Edmonton, Alberta \\
\> Canada T6G 2G1 \\[\medskipamount]
{\it E-mail\/}: \> {\tt runde@math.ualberta.ca} \\
\> {\tt vrunde@ualberta.ca} \\[\medskipamount]
{\it URL\/}: \> {\tt http://www.math.ualberta.ca/$^\sim$runde/runde.html}
\end{tabbing} 
\end{document}

%% file: mathdefs.tex
\usepackage{amsmath,amssymb,amscd}
\newcommand{\pf}[1]{\trivlist \item[\hskip\labelsep\it #1\ ]}
\newcommand{\varpf}[1]{\trivlist \item[\hskip\labelsep\sc #1:]}
\newcommand{\qedbox}{$\rlap{$\sqcap$}\sqcup$}
\newcommand{\qed}{\qquad \qedbox \endtrivlist}
\newcommand{\varqed}{\hfill \rule{0.6em}{0.6em} \endtrivlist}

\newenvironment{items}{
  \begin{enumerate} 
                    
  }{\end{enumerate}}